# A New Algorithm for Multicommodity Flow

Dhananjay P. Mehendale
Sir Parashurambhau College, Tilak Road, Pune 411030, India

## Abstract

We propose a new algorithm to obtain max flow for the multicommodity flow. This algorithm utilizes the max-flow min-cut theorem and the well known labeling algorithm due to Ford and Fulkerson [1]. We proceed as follows: We select one source/sink pair among the n distinguished source/sink pairs at a time and treat the given multicommodity network as a single commodity network for such chosen source/sink pair. Then applying standard labeling algorithm, separately for each sink/source pair, the feasible flow which is max flow and the corresponding minimum cut corresponding to each source/sink pair is obtained. A record is made of these cuts and the paths flowing through the edges of these cuts. This record is then utilized to develop our algorithm to obtain max flow for multicommodity flow. In this paper we have pinpointed the difficulty behind not getting a max flow min cut type theorem for multicommodity flow and found out a remedy.

**1. Introduction:** Multicommodity flow problem arises in a wide variety of important applications. Many communications, logistics, manufacturing, and transportation problems can be formulated as large multicommodity flow problems. The max-flow min-cut theorem due to Ford and Fulkerson [1] is applicable to transport networks containing one source and one sink, represented by simple connected weighted digraphs. However, if there are several sources $s_1, s_2, \cdots, s_k$ and several sinks $t_1, t_2, \cdots, t_k$ forming a directed network and if there exists a restriction that the flow from the specified source, $s_i$, should be sent to the specified sink, $t_i$, then this problem is called the problem of multicommodity flow. There is no result similar to max-flow min-cut theorem for multicommodity flow problem in general. This paper aims at pinpointing the reason behind not getting a max flow min cut type theorem for multicommodity flow and further a way to surmount this difficulty.

In multicommodity flow problems several distinct commodities are flowing simultaneously through the given transport network formed by directed edges connecting vertices of the network and each commodity has its own source and own sink. All flows share the edge capacities and as in the case of single commodity, having single source and single sink, the sum of all flows through an edge must not exceed the capacity of the edge. Also, the conservation conditions, as in the case of single commodity case, are observed at each intermediate node.

**2. Preliminaries:** For network with one source and one sink there exists very well known algorithm due to Ford and Fulkerson to obtain max flow [1]. This so called labeling algorithm starts with a feasible flow as an input and proceeds to find an $f$- augmenting path $P$ with positive leeway $\varepsilon$ and then increasing the flow by $\varepsilon$ along forward edges of $P$ and decreasing the flow by $\varepsilon$ along backward edges of $P$ produces a feasible flow $f'$ with value higher by $\varepsilon$, i.e. value($f'$) = value($f$) + $\varepsilon$. If this labeling algorithm terminates at some iteration then it terminates (only) when it does not find any $f$- augmenting path $P$. This termination (in finitely many iterations) always takes place, i.e., starting with node set $S = \{s\}$, where element $s$ is source node, we reach at nodes, which get labels by the algorithm, and from these labeled nodes we cannot further label any nodes and the labeling algorithm thus terminates producing a cut $[S, S']$ where $S'$ is nonempty complement of $S$ in the entire node set. In short, the labeling algorithm always terminates after finitely many iterations by landing at a cut (or cutset) with same value as this flow. By weak duality theorem, namely, if $f$ is a feasible flow and $[S, T]$ is a source/sink cut, then value($f$) $\leq$ cap($[S, T]$), where cap stands for capacity, the cut is a minimum cut and the flow is a maximum flow, i.e. both are optimal (See 4.3, page 158, [2]).

**3. Basic Preparation for Algorithm:** The input for this algorithm is the given multicommodity network, a digraph with a nonnegative capacity, $c(e)$, on each (directed) edge $e$ and there are $n$ distinguished source vertices (or nodes) $s_1, s_2, \cdots, s_n$ and corresponding $n$ distinguished sink vertices (or nodes) $t_1, t_2, \cdots, t_n$ with a restriction that the flow from the specified source, $s_i$, should be sent to the specified sink, $t_i$. We begin with describing certain steps to be carried out on given multicommodity network and we call this processing "Initialization".

**Step 1:** We begin with treating this multicommodity network as a single commodity network with a single source vertex (or node) $s_1$ and single sink vertex (or node) $t_1$, and treat all other nodes as ordinary nodes as in the sense of usual one source and one sink network.

**Step 2:** We apply the algorithm of Ford and Fulkerson to the network, which is now treated as single source/sink network arrived at step 1, and (in finitely many iterations) find max flow and min cut.

**Remark 3.1:** We want to repeat these steps 1, 2 with all other source/sink pairs ($s_i, t_i$), $i = 2, 3, \ldots, n$ and at each carrying out of these steps 1, 2 for each source/sink pair we want to keep record of the max flow we get, the min cut we get, the paths that are flowing this max flow we get, and finally the edges through which these paths pass. We actually want to use in our algorithm this information mentioned here that is associated with each source/sink pair. So, we proceed as follows:

**Step 3:** At the end of step 2 we make a record of distinct paths which are starting at source vertex $s_1$, then passing through some edge of min cut, and terminating in sink vertex $t_1$. We assign a distinct color with each of these distinct paths and further the color which is assigned to a path is also assigned to each edge on that path. So, let $\{P_1^1, P_2^1, \ldots, P_{r_1}^1\}$ be the paths starting at source vertex $s_1$, then passing respectively through edges $\{e_1^1, e_2^1, \cdots, e_{r_1}^1\}$ = min cut, and terminating in sink vertex $t_1$. Also, let $\{e_{j1}^1, e_{j2}^1, \cdots, e_{ju_1}^1\}$ are the edges on paths $P_j^1$, $j = 1, 2, \cdots, r_1$. We associate distinct colors $\{c_1^1, c_2^1, \cdots, c_{r_1}^1\}$ with paths $\{P_1^1, P_2^1, \ldots, P_{r_1}^1\}$ respectively and also associate same colors with each edge on a path, thus, color $c_j^1$ is associated with each edge among the edges $\{e_{j1}^1, e_{j2}^1, \cdots, e_{ju_1}^1\}$ on path $P_j^1$.

**Step 4:** For keeping this coloring record in a systematic and convenient way, we then proceed, for the given multicommodity network *N*, with construction of a bitableau, called the Edge-Color-Bitableau, *ECBT(N)*, to keep the record of colors we have associated with edges in step 3, in a systematic way for using it later as follows:

Initially, we have (before assignment of any color) a left tableau consisting of a column vector representing edge labels of edges in N and a blank right tableau:

$$ECBT(N) = \begin{pmatrix} e_1 & \cdots \\ e_2 & \cdots \\ \vdots & \\ e_j & \cdots \\ \vdots & \\ e_p & \cdots \end{pmatrix}$$

After step 3, we have assigned colors with certain edges that appeared on certain path among the paths. We write down the colors assigned to an edge in the same row in the right tableau in which that edge label appears in the left tableau. We carry out this procedure for every edge on every path $P_j^1$, $j = 1, 2, \cdots, r_1$ and thus partially build the initially empty right tableau with the corresponding color labels associated with the edges involved.

**Step 5:** We now go back to step 1. We replace there source/sink pair $(s_1, t_1)$ by new source/sink pair $(s_2, t_2)$ and proceed again with steps 1, 2, 3, 4 leading to new paths $\{P_1^2, P_2^2, \ldots, P_{r_2}^2\}$ with $\{e_1^2, e_2^2, \cdots, e_{r_2}^2\}$ = min cut, and $\{e_{j1}^2, e_{j2}^2, \cdots, e_{ju_2}^2\}$ are the edges on paths $P_j^2$, $j = 1, 2, \cdots, r_2$. We associate distinct colors $\{c_1^2, c_2^2, \cdots, c_{r_2}^2\}$ with paths $\{P_1^2, P_2^2, \ldots, P_{r_2}^2\}$ respectively and also associate same colors with each edge on a path, thus, color $c_j^2$ is associated with each edge among the edges $\{e_{j1}^2, e_{j2}^2, \cdots, e_{ju_2}^2\}$ on path $P_j^2$. We write down the colors assigned to an edge in the same row in the right tableau in which that edge label appears in the left tableau. We carry out this procedure for every edge on every path $P_j^2$, $j = 1, 2, \cdots, r_2$ and thus continue with further partially building of the right tableau with the corresponding color labels associated with the edges involved this time.

**Step 6:** We now proceed with repeating step 5 for every other source/sink pair $(s_i, t_i)$, $i = 3, 4, \ldots, n$, and thus, complete the construction of *ECBT(N)*.

**Step 7:** We now proceed, for the given multicommodity network *N*, with construction a bivector called the Edge-Capacity-Bivector, *ECBV(N)*, one more bitableau called the Path-Recorder-Bitableau, *PRBT(N)*, one more bivector called Minimum-Edge-Capacity-Bivector, *MECBV(N)*, and one more bivector called Path-Cardinality-Bivector, *PCBV(N)*, for using it later to make decision about choosing paths to maximize the flow as follows:

$$ECBV(N) = \begin{pmatrix} e_1 & Cap(e_1) \\ e_2 & Cap(e_2) \\ \vdots & \\ e_j & Cap(e_j) \\ \vdots & \\ e_p & Cap(e_p) \end{pmatrix}$$

,

$$PRBT(N) = \begin{pmatrix} P_1^1 & e_{11}^1(w_{11}^1) & \cdots & e_{1j}^1(w_{1j}^1) & \cdots & e_{1u_1}^1(w_{1u_1}^1) \\ P_2^1 & e_{21}^1(w_{21}^1) & \cdots & e_{2j}^1(w_{2j}^2) & \cdots & e_{2u_2}^1(w_{2u_2}^1) \\ \vdots & & & & & \\ P_j^i & e_{j1}^i(w_{j1}^i) & \cdots & e_{jj}^i(w_{jj}^i) & \cdots & e_{ju_j}^i(w_{ju_j}^i) \\ \vdots & & & & & \\ P_{r_n}^n & e_{r_n 1}^n(w_{r_n 1}^n) & \cdots & e_{r_n j}^n(w_{r_n j}^n) & \cdots & e_{r_n u_{r_n}}^n(w_{r_n u_{r_n}}^n) \end{pmatrix}$$

$$MECBV(N) = \begin{pmatrix} P_1^1 & & Cap(e_1^1) \\ P_2^1 & & Cap(e_2^1) \\ \vdots & & \\ P_j^i & & Cap(e_j^i) \\ \vdots & & \\ P_{r_n}^n & & Cap(e_{r_n}^n) \end{pmatrix}$$

and,

$$PCBV(N) = \begin{pmatrix} P_1^1 & & Card(E_1^1) \\ P_2^1 & & Card(E_2^1) \\ \vdots & & \\ P_j^i & & Card(E_j^i) \\ \vdots & & \\ P_{r_n}^n & & Card(E_{r_n}^n) \end{pmatrix}$$

a) The left vector of *ECBV(N)* is a column vector consisting of edge labels of the edges of network *N*.
b) The right vector of *ECBV(N)* is a column vector consisting of capacities of the corresponding edge labels written in the left vector in the same row.
c) The left tableau of *PRBT(N)* is a column vector consisting of path labels of the paths, $P_j^i$, $j = 1, 2, \cdots, r_i$ and $i = 3, 4, \ldots, n$, obtained in succession for all source/sink pair $(s_i, t_i)$ of network *N*.
d) In the right tableau of *PRBT(N)* the vector of edge labels, along with their capacities in bracket in front of them, are written in the same order as they have occurred on these paths.
e) The left vector of *MECBV(N)* is same as left tableau of *PRBT(N)*.
f) The right vector contains vector of minimum capacities, i.e. with minimum capacity on the paths whose path labels are written in the

same row in the left vector. It is clear that initially the value will be the value of the capacity of the edge belonging to minimum cut on that on the path.

g) The left vector of *PCBV*(*N*) is same as the left tableau of *PRBT*(N).
h) The right vector of *PCBV*(*N*) is a column vector consisting of count of the distinct colors associated with the edges on this path. Note that these colors associated with edges through their appearance on different paths is depicted in *ECBT*(*N*). Some more explanation is as follows: Let the directed edges present on path $P_j^i$ be $\{e_{j1}^i, e_{j2}^i, \cdots, e_{ju_i}^i\}$. Now look in those rows of right tableau of *ECBT*(*N*) where there are edge labels $e_{j1}^i, e_{j2}^i, \cdots, e_{ju_i}^i$ present in the left tableau of these rows and collect labels of all distinct (without repetition) colors in a set, say, $E_j^i$. Let *Card* ($E_j^i$) denotes the cardinality of this set of distinct color labels associated together with the edges on the path $P_j^i$ and so these distinct color labels together (collected in $E_j^i$) are essentially associated with the path $P_j^i$.

The basic preparation for the algorithm is now complete and we are now ready for the discussion of the actual steps of the algorithm.

**4. Algorithm for Multicommodity Flow:** Let us note down first some important observations:

1. We have associated a unique color with each path and also that same color with each edge on this path.
2. Same edge can appear on more than one path and therefore more than one color can get associated with the same edge.
3. When we select some path, $P_j^i$, to pass the commodity, say, $M_i$, of quantity, say, $Cap(e_j^i)$, where $e_j^i$ is the edge on the path $P_j^i$ belonging to min cut obtained for pair $(s_i, t_i)$, then we cannot use the same edge to pass some other commodity $M_k$ if the same edge appears as edge $e_l^k$ on other path $P_l^k$ on the min cut obtained for pair $(s_k, t_k)$.
4. The maximum possible flow on a path is equal to the minimum of the capacities of the edges belonging to that path, i.e. if the directed edges present on path $P_j^i$ are $\{e_{j1}^i, e_{j2}^i, \cdots, e_{ju_i}^i\}$ and if edge $e_j^i \in \{e_{j1}^i, e_{j2}^i, \cdots, e_{ju_i}^i\}$ has minimum capacity, or, is the one which belongs to corresponding min cut, then max flow that is possible on this path = $Cap(e_j^i)$.

5. When the quantity of commodity with value = $Cap(e_j^i)$ corresponding to path $P_j^i$ is transported the "free (or, usable) capacity" for the other edges on this path becomes $|Cap(e_{jm}^i) - Cap(e_j^i)|$, $m = 1, \cdots, u_i$.

6. If more than one color gets associated with an edge due to its appearance on more than one path then we can use that edge with its full capacity only for some one and one color among these associated colors, i.e. for flowing only some one commodity.

7. If more than one color gets associated with an edge due to its appearance on more than one path then we can use that edge with partial capacity for multiple colors among these associated colors, i.e. for transporting more than one commodity, by observing the capacity constraint for that edge.

8. If $Card(E_j^i) = 1$ for all $i$ and $j$ in the right vector of $PCBV(N)$, i.e. if the right vector of $PCBV(N)$ is made up of ones, then it meant that there is no sharing of edges while flowing the commodities through the network and we can achieve max flow for each commodity separately (as per max flow min cut theorem). Therefore, the total flow, equal to sum of individual max flows, will be the max flow for multicommodity network, $N$.

9. In general, $Card(E_j^i) > 1$ for many pairs of $i$ and $j$. This meant same edge appears on many paths. This essentially implies that many paths try to share this same edge for transporting their associated commodities for achieving their individual max flow through max flow min cut algorithm.

10. "Max flow" possible for multicommodity flow is equal to $f$, where,
$$f = \sum_i S(C_i) + (-1)^1 \sum_{i<j} S(C_i \cap C_j) + \cdots + (-1)^{(n-1)} S(C_1 \cap C_2 \cap \cdots \cap C_n)$$
and where, $S(C_i)$ = sum of capacities of edges of min cut $C_i$.

11. The question that naturally arises is: How to achieve this maximum possible flow? Which path among the paths that passes through a certain edge, common to them, should be chosen to flow its associated commodity and thus use the capacity of this edge which is member of these paths? The simple and natural answer is: Choose that path (or that color), any one among the paths that pass through the edge which is common to them, to flow the commodity, which in this process eliminates minimum number of paths for further transport.

**Steps of the Algorithm:**
1. Carry out the process of "Initialization" on the given mulicommodity network, i.e. carry out the process described under the heading "Initialization" in the steps given there.
2. From *PCBV(N)* find out the path labels $P_j^i$ in the left vector for which the values of $Card(E_j^i)$ in right vector are unity and transport the commodities on these paths and mark these paths as used for transport. Make the record of total transport achieved in this process.
3. From *PCBV(N)*, among the unmarked rows after the first step, find out the path labels $P_j^i$ in the left vector of *PCBV(N)* for which the values of $Card(E_j^i)$ in right vector is minimum. Among such rows select any row with path label say, $P_j^i$, to transport the commodity associated with thus chosen path with value equal to value of minimum capacity edge on that path, say, $Cap(e_j^i)$, and mark this path as used for transport. Modify the record of total transport achieved in this process. Subtract this value of minimum capacity from the original capacity of every edge that has occurred on this path. Modify *ECBV(N)*, *PRBT(N)*, *MECBV(N)*, *PCBV(N)*. If this modification produces any paths containing some edge with zero (left out) capacity then discard those paths from further selections, by marking such paths as discarded.
4. Continue step 3 till every path is either is marked as selected or discarded.
5. Record total transport achieved in this process, carried out till end. It will be the "max flow" for the given multicommodity network under consideration.

**Example:** Consider following multicommodity network with two commodities, as depicted in FIG. 2.1. The capacities are written near the edges.

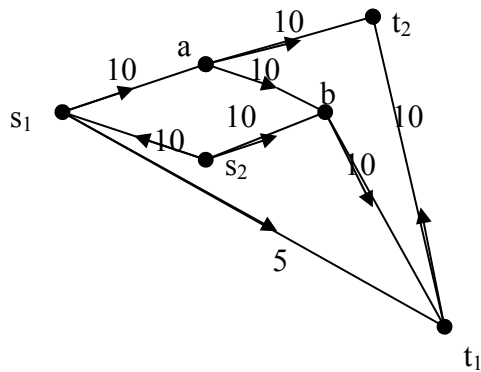

FIG. 2.1

This network has two source/sink pairs $(s_1, t_1)$ and $(s_2, t_2)$.
As per the "Initialization" procedure, the paths forming *PRBT(N)* are:

$$P_{11}^1 : s_1 \xrightarrow{5} t_1 : \text{Violet}$$

$$P_{12}^1 : s_1 \xrightarrow{10} a \xrightarrow{10} b \xrightarrow{10} t_1 : \text{Red}$$

$$P_{11}^2 : s_2 \xrightarrow{10} s_1 \xrightarrow{10} a \xrightarrow{10} t_2 : \text{Green}$$

$$P_{12}^2 : s_2 \to b \to t_1 \to t_2 : \text{Yellow}$$

where on the top of every edge the capacity is depicted. Also, we have written the colors applied to paths initially on the right side of the path. Using these paths and directed edges and their associated colors and capacities, We can build *ECBT(N)* and using it we can further build *ECBV(N)*, *PRBT(N)*, *MECBV(N)*, *PCBV(N)*.
When we proceed as per algorithm we see that finally (from construction of *PCBV(N)*):

- Count of colors associated with $P_{11}^1 = 1$

- Count of colors associated with $P_{12}^1 = 3$
- Count of colors associated with $P_{11}^2 = 2$
- Count of colors associated with $P_{12}^2 = 2$

Therefore, by proceeding as per algorithm:

1) We choose first path $P_{11}^1$
2) We then choose path $P_{11}^2$
3) We then choose path $P_{12}^2$
4) We discard path $P_{12}^1$

And we achieve max flow = 5 units + 10 units + 10 units = 25 units. Out of these units, 5 units are of the first commodity and 20 units are of the second commodity.

**Theorem:** Algorithm for multicommodity flow produces max flow.

**Proof:** The algorithm begins with applying max flow min cut theorem individually and separately to each source/sink pair of the given multicommodity flow and produces separate record of paths, min cuts, etc. producing max flow for individual commodities. When there is no sharing of edges and thus no intersection of these paths then the sum of individual max flows will be the max flow for multicommodity network.
When there is sharing of edges by different paths the algorithm finds out the optimal way of selecting paths among these paths for flowing the commodities so that minimum numbers of paths get eliminated from the totality of paths and thus causes minimum possible reduction in the sum of individual max flows, and thus, produces the desired max flow for multicommodity network under consideration.

□